\documentclass[11pt,twoside]{article}

\usepackage{amsbsy,amsfonts,amsmath,amssymb,eucal}
\usepackage[all]{xy}

\def\v#1{\par \vspace{#1mm} \par}
\def\itemn#1{\item[\hspace{0.6mm} {\rm (#1)}]}

\def\qedend{\hfill $\square$}
\def\bin#1#2{\mbox{\scriptsize $\big\{\!\!\!\!\begin{array}{c} #1 \\ #2
    \end{array}\!\!\!\!\big\}$\normalsize}}
\renewcommand{\geq}{\geqslant}
\renewcommand{\leq}{\leqslant}
\def\sm{{\rm sm}}
\def\zmod#1{\mathbb{Z}/#1\mathbb{Z}}
\def\et{{\rm \acute{e}t}}
\def\longto{\longrightarrow}
\def\isomto{\stackrel{\sim}{\longto}}

\renewcommand{\tilde}{\widetilde}
\renewcommand{\hat}{\widehat}

\newtheorem{definition}{Definition}[subsection]
\newenvironment{defi}{\begin{definition} \rm}{\end{definition}}
\newtheorem{prop}[definition]{Proposition}
\newtheorem{lemm}[definition]{Lemma}
\newtheorem{coro}[definition]{Corollary}
\newtheorem{theo}[definition]{Theorem}
\newtheorem{construction}[definition]{Construction}
\newenvironment{cons}{\begin{construction} \rm}{\end{construction}}
\newtheorem{remark}[definition]{Remark}

\newtheorem{example}[definition]{Example}
\newenvironment{exam}{\begin{example} \rm}{\end{example}}
\newtheorem{nothing}[definition]{$\!\!$}
\newenvironment{noth}{\begin{nothing} \rm}{\end{nothing}}
\newenvironment{proo}{{\flushleft \bf Proof :}}{\hfill $\square$ \vspace{5mm}}

\DeclareMathOperator{\Aut}{Aut}
\DeclareMathOperator{\Id}{id}
\DeclareMathOperator{\Spec}{Spec}

\def\clF{{\cal F}} \def\clG{{\cal G}} \def\clH{{\cal H}} \def\clI{{\cal I}}
\def\clK{{\cal K}}   \def\clN{{\cal N}}
\def\clO{{\cal O}}
   
   \def\clZ{{\cal Z}}

\def\bbA{{\mathbb A}} \def\bbQ{{\mathbb Q}} \def\bbZ{{\mathbb Z}}

\begin{document}

\begin{center}
{\bf \Large Effective model of a finite group action}
\v{10}{\bf Matthieu Romagny}
\v{2}
{\em January 20, 2006}
\end{center}

\v{7}

{\def\thefootnote{\relax}
\footnote{ \hspace{-6.8mm}
}
}

Let $R$ be a discrete valuation ring with fraction field $K$,
uniformizer $\pi$ and residue field $k$ of characteristic $p>0$. Let
$X_K$ be a flat $K$-scheme of finite type endowed with a faithful
action of a finite group scheme $G$. Given an $R$-model $X$ to which
the action extends, it may happen that the reduced action on the
special fibre acquires a kernel, especially if $p$ divides
$|G|$. Typically, this happens if $K$ has characteristic $0$ and
$X_K=A_K$ is an abelian variety of dimension $g$ with rational
$p$-torsion, $G=(\zmod{p})^{2g}$ acts by translations, and $X=A$ is
the N{\'e}ron model. However, in this example we feel that there is a
better group acting~: namely the $p$-torsion subgroup $A[p]$ acts
faithfully on both fibres of $A$. In general, under the natural
assumptions below, we will prove that there always exists a group
$\clG$ faithful on the special fibre, coming with a dominant map $G\to
\clG$ which is an isomorphism on the generic fibre. In our setting, we
shall actually consider that a model $X$ with an action of $G$ is
given from the start, and call $\clG$ the {\em effective model} for
the action.

\medskip

\noindent {\bf Theorem~:} {\em Let $G$ be a finite flat $R$-group
  scheme. Let $X$ be a flat $R$-scheme of finite type with an action
  of $G$. Assume that $X$ is covered by $G$-stable open affines $U_i$
  with function ring separated for the $\pi$-adic topology, such that
  $G_K$ acts faithfully on the generic fibre $U_{i,K}$. Then if $X$
  has reduced special fibre, there exists an effective model for the
  action of $G$.}

\medskip

This is corollary \ref{coroMT} below.
We now briefly put our result in perspective by describing its
original motivation and some of its corollaries related to other
current work.

Our main motivation comes from Galois covers of curves. Assume that
$K=\bbQ_p$ is the field of $p$-adic numbers. Then there is a nice
smooth (nonconnected) proper stack classifying admissible $G$-Galois
covers of stable curves, with fixed ramification invariants. At the
moment, the question of understanding its reduction at $p$ seems wide
open~; one explanation is the following. A natural $R$-model for a
smooth curve $X_K$ is the stable model. If $G$ acts on $X$, the above
phenomenon where the group action degenerates shows up locally on
components of the special fibre $X_k$. This perturbates the usual
local-global principle saying that deformations are localized at
singular points and ramification points... Thus the theorem is an
attempt to remedy this pathology.

There are other useful notions of effective models. In
fact let $X\to Y=X/G$ be a finite cover of stable $R$-curves. If
$p^2{\not |}\: |G|$ and the $p$-Sylow is normal, then Dan Abramovich
has shown in~\cite{Ab} that there exists a finite flat group scheme
$\clG'\to Y_\sm$ on the smooth locus of the base, acting on
$X_\sm$ in such a way that $X_\sm\to Y_\sm$ is a torsor away from some
relative divisor. Abramovich calls it {\em Raynaud's group scheme} and
shows how one may twist $Y$, in the sense of twisted curves, so as to
extend $\clG$ to the whole curve. Note that since $\Aut_R(X)$ is
unramified, the effective model obtained from our theorem is $\clG=G$,
so $\clG'$ is not just the pullback to $Y_\sm$ of ours. However, one
salient feature of our construction is that no properness of $X$ is
required, so that we may perform it locally on $Y$ and glue to recover
the model $\clG'$. We obtain some additional information that can not
be reached from \cite{Ab}~: for example we prove that Raynaud's group
scheme is constant on irreducible components of $Y_k$, and also that
it extends to nodes $p\in Y_k$ lying on one single component (without
the need of a twisting), see corollary~\ref{precisions_on_RGS}.

The main theorem occupies section~1 of the paper. Finally we remark
that the construction of Raynaud's group scheme in \cite{Ab} can be
carried out without assumptions on $|G|$, as well as ours. The trouble
is that in general one does not get a torsor structure on the special
fibre, as is known to the experts (see Sa{\"\i}di \cite{Sa1} for an
example). Our effective model gives an explanation for this,
illustrated by several examples exposed in section~2.

\medskip

{\em Acknowledgements}. I wish to thank Jos{\'e} Bertin for raising this
question in the context of my thesis and for reading a first
version of this note. Laurent Moret-Bailly suggested the lines of a
strategy of proof and is also warmly acknowledged. Finally I am
thankful to Dan Abramovich for inviting me to Brown University in
march 2005, and for several remarks and counterexamples leading to
more accurate statements.

\section{Models for finite group actions}

Throughout we fix a discrete valuation ring $R$ with fraction field $K$,
uniformizer $\pi$ and residue field $k$ of characteristic $p>0$. For
$G$ a group scheme over $R$, we denote by $R[G]$ or simply $RG$ its
ring of functions.

\subsection{Definitions, basic properties} \label{defiRappels}

We first recall some well-known properties of the scheme-theoretic image for
morphisms over a discrete valuation ring $R$, and we fix some terminology.

\begin{trivlist}
\itemn{i} Let $f\colon W\to X$ be a morphism of $R$-schemes such that
$f_*\clO_W$ is quasi-coherent. Then there exists a smallest closed subscheme
$X'\subset X$ such that $f$ factors through $X'$. We call it the
schematic image of $f$. If it is equal to $X$ we say that $f$ is dominant.
\itemn{ii} If $W=Z_K$ is a closed subscheme of the generic fibre of $X$
and $f$ is the canonical immersion, then the schematic image is called the
schematic closure of $Z_K$ in $X$. It is the unique closed subscheme
$Z\subset X$ which is flat over $R$ and satisfies $Z\otimes K=Z_K$.
\itemn{iii} Let $X$ be a scheme over $R$. A family of closed subschemes
$Z_\alpha\subset X$ with ideal sheaves $\clI_\alpha$ is called schematically
dense in $X$ if the morphism $\amalg Z_\alpha\to X$ is dominant. In
other words this means that $\cap\, \clI_\alpha=0$ in $\clO_X$. Any
other family $Z'_\beta\subset X$ such that for any $\alpha$ there
exists one $\beta$ with $Z_\alpha\subset Z'_\beta$, is again
schematically dense in $X$.
\itemn{iv} Let $G$ be a finite, flat group scheme over $R$ and let $X$ be an
$R$-scheme. Let $\mu\colon G\times X\to X$ be an action. Then $\mu$ is a
finite, flat morphism of finite presentation. If $Z\subset X$ is a closed
subscheme we denote by $G.Z$ the schematic image of $G\times Z$ under $\mu$. If
$Z$ is finite (resp. flat) over $R$, then $G.Z$ also is.
\end{trivlist}

We now recall from \cite{Ra2}, \S~2 the relation of \emph{domination}
between models of a group~:

\begin{defi}
\begin{trivlist}
\itemn{i} Let $G_1$ and $G_2$ be finite flat group schemes over $R$
with an isomorphism $u_K\colon G_{1,K}\to G_{2,K}$. We say that $G_1$
\emph{dominates} $G_2$ and we write $G_1\geq G_2$, if we are given an
$R$-morphism $u\colon G_1\to G_2$ which restricts to $u_K$ on the
generic fibre. If moreover $G_1$ and $G_2$ act on $X$, we say that
$G_1$ dominates $G_2$ \emph{compatibly} (with the actions)
if $\mu_1=\mu_2\circ (u\times\Id)$.
\itemn{ii} Let $G$ be a finite flat group scheme over $R$. Let $X$ a flat
scheme over $R$. Let $\mu\colon G\times X\to X$ be an action, faithful on the
generic fibre. An \emph{effective model} for $\mu$ is a finite flat $R$-group
scheme $\clG$ acting on $X$, dominated by $G$ compatibly, such that $\clG$
acts universally faithfully on $X$, that is to say, faithfully both on
the generic and the special fibre.
\end{trivlist}
\end{defi}

\begin{lemm}
An effective model is unique up to unique isomorphism, if it exists.
\end{lemm}

\begin{proo}
Giving the action is the same as giving a morphism of fppf sheaves
$f\colon G\to \Aut_R(X)$. The requirements of universal effectivity and of
flatness say that if it exists, $\clG$ must be the scheme-theoretic closure of
$G_K$ in the fppf sheaf $\Aut_R(X)$, as defined in \cite{Ra1}, \S~3.
\end{proo}

We recall that an action $\mu\colon G\times X\to X$ is called \emph{admissible}
if $X$ can be covered by $G$-stable open affine subschemes.

\begin{prop} \label{propUnicite}
Let $G$ be a finite flat group scheme over $R$. Let $X$ a flat scheme over
$R$ and $\mu\colon G\times X\to X$ an admissible action, faithful on the generic
fibre. Assume that there exists an effective model $G\to \clG$. Then,
\begin{trivlist}
\itemn{i} If $H$ is a finite flat subgroup of $G$, the restriction of the action
to $H$ has an effective model $\clH$ which is the schematic image of $H$ in $\clG$.
If $H$ is normal in $G$, then $\clH$ is normal in $\clG$.
\itemn{ii} The identity of $X$ induces an isomorphism $X/G\simeq X/\clG$.
\itemn{iii} Assume that $G$ is {\'e}tale and note $p={\rm char}(k)$. Let
$N\lhd G$ be the (unique) subgroup of $G$ such that $N_k$ is the kernel
of the action on $X_k$. Then, its effective model $\clN$ is a connected
$p$-group.
\itemn{iv} Assume that there is an open subset $U\subset X$ which is schematically dense
in any fibre of $X\to\Spec(R)$, such that $\clG$ acts freely on $U$. Then
for any closed normal subgroup $H\lhd G$, the effective model
of $G/H$ acting on $X/H$ is $\clG/\clH$.
\itemn{v} Under assumptions (iii)+(iv) the model $\clG$ has a connected-{\'e}tale
sequence
$$1\to \clN \to \clG\to G/N \to 1$$
\end{trivlist}
\end{prop}

\begin{proo}
(i) is clear.

\smallskip

\noindent (ii) Locally, $X=\Spec(A)$ and $X/G$ is defined as the spectrum
of the ring $A^G=\{\,a\in A,\,\mu_G^\sharp(a)=1\otimes a\,\}$ where
$\mu_G^\sharp$ is the coaction. Now $\mu_G^\sharp$ factors through the
coaction $\mu_\clG^\sharp$ corresponding to the action of $\clG$ :
$$
A\to R\clG\otimes A\hookrightarrow RG\otimes A
$$
Therefore, $A^\clG=\{\,a\in A,\,\mu_\clG^\sharp(a)=1\otimes a\,\}=A^G$.

\smallskip

\noindent (iii)
Since the composition $N_k\to \clN_k\hookrightarrow \Aut_k(X_k)$ is trivial
(as a morphism of sheaves), the morphism $N_k\to \clN_k$ also is. Moreover,
$N\to \clN$ is dominant and closed hence surjective.
Hence $\clN_k$ is infinitesimal so $\clN$ is a
$p$-group. Let us show that it is connected. We may and do assume that
$R$ is complete. Then $\clN$ has a connected-{\'e}tale sequence,
we denote the {\'e}tale quotient by $\clN_\et$.
The composite map $t:N\to \clN\to \clN_\et$ is trivial on the
special fibre. Moreover, $t$ is determined by its restriction to the special
fibre because it is a morphism between {\'e}tale schemes. So it is globally
trivial. As $t$ is dominant we get $\clN_\et=1$ thus $\clN$
is connected.

\smallskip

\noindent (iv) Clearly $H$ acts admissibly, and $X/H\isomto X/\clH$
by (ii). We just have to show that
$\clG/\clH$ acts universally faithfully on $X/\clH$. Under the assumption
of existence of $U\subset X$, one checks easily that this follows from the fact
that $\clG$ acts universally faithfully on $X$.

\smallskip

\noindent (v) Apply (iv) to $H=N$.
\end{proo}

\subsection{Existence}

In this section, our aim is to prove that an effective model exists
when the scheme $X$ is of finite type over $R$ with
reduced special fibre (corollary \ref{coroMT} below).
For the moment let $G$ be a finite flat group scheme over $R$ and
let $X$ be a flat $R$-scheme on which $G$ acts (no reducedness
assumption). As always we assume that the action $\mu$ is faithful on
the generic fibre.

\begin{lemm} \label{lemmProper}
If $X$ is proper, there exists an effective model for the action of $G$.
\end{lemm}

\begin{proo}
In this case, Artin has shown that the fpqc sheaf of automorphisms
$\Aut_R(X)$ is an algebraic space locally of finite presentation.
Moreover, $X$ being flat over $R$, it is separated. It follows
that $\Aut_R(X)$ is actually a scheme (\cite{An}, chapter IV). By
assumption we have a morphism $u:G\to \Aut_R(X)$ and $G_K$ is a closed
subscheme of $\Aut_K(X_K)$.
We call $\clG$ the schematic image of $u$. It is finite and flat over $R$. By
its definition, $\clG$ acts universally faithfully on $X$.
\end{proo}

We denote by $\clF$ the family of all closed subschemes $Z\subset X$
which are finite and flat over $R$. 

\begin{theo} \label{theoZ}
Let $X$ be a flat $R$-scheme and $\mu\colon G\times X\to X$ an action.
Assume that $X$ is covered by $G$-stable open affines $U_i$ with function
ring separated for the $\pi$-adic topology, such that $G$ acts faithfully
on the generic fibre $U_{i,K}$.
Assume that $\clF_k$ (the family of all $Z_k$) is schematically dense
in $X_k$. Then there exists an effective model for the action of $G$.
\end{theo}

\begin{proo}
The proof is divided into three steps~:
\begin{trivlist}
\itemn{1} The subfamily $\clF^*$ of all $Z$ such that $G$ stabilizes
$Z$ and acts faithfully on $Z_K$ is again schematically dense in the special
fibre $X_k$.
\end{trivlist}
The point here is that, by lemma \ref{lemmProper}, for any $Z\in \clF^*$ there
is an effective model $G\to \clG_Z$.
\begin{trivlist}
\itemn{2} There exists $\clZ\in \clF^*$ such that $\clG:=\clG_\clZ$ acts
on all subschemes $Z\in \clF^*$.
\itemn{3} The group $\clG$ is the effective model of the action of $G$.
\end{trivlist}

Let us prove (1).
If $Z\in \clF$ then the orbit $G.Z$ (see \ref{defiRappels}(iv))
is again in $\clF$. Also, $\clF$ is stable under union of
subschemes, in the sense of intersection of defining ideals. Moreover,
$\clF^*$ contains at least one element $Z_0$ : consider a finite
$K$-subscheme of $X_K$, on which $G_K$ acts faithfully, and look at
its schematic closure in $X$. Therefore $\clF^*$ contains all $G.Z\cup Z_0$,
for $Z\in \clF$, so a fortiori $\clF^*_k$ is schematically dense in $X_k$
(see \ref{defiRappels}(iii)).

We now prove (2).
As already said, $\clF$ is a filtering set with the order given by
inclusion. If $Z_1\subset Z_2$ then $\clG_{Z_2}$ stabilizes $Z_1$ (because
this is true on the generic fibre). Hence the effective model for
$\clG_{Z_2}$ acting on $Z_1$ can be identified with $\clG_{Z_1}$,
and there is a map
$$
\clG_{Z_2}\to \clG_{Z_1}
$$
This gives inclusions of structure rings $R\clG_{Z_1}\subset R\clG_{Z_2}
\subset RG$. As $RG$ is a finite $R$-module there exists $\clZ\in
\clF^*$ such that $R\clG_\clZ$ contains all the other $R\clG_Z$.
Therefore $\clG_{\clZ}$ acts on all $Z$'s.

It only remains to prove (3).
The point is to see that $\clG$ acts on $X$. In view of the
assumptions of the theorem, in order to construct this action we can
assume that $X$ is affine and $G$-stable, equal to $\Spec(A)$ with $A$
separated for the $\pi$-adic topology. For $Z\in \clF^*$ denote by
$I_Z$ its ideal in $A$. Consider the flat $R$-algebra
$$
\hat{A}:=\lim_{\longleftarrow \atop Z\in \clF^*} A/I_Z
$$
The kernel of the natural map $A\to \hat{A}$ is the intersection of all $I_Z$.
Under reduction modulo $\pi$, this intersection maps inside $\cap\,(I_Z\otimes k)$.
By step (1) this is zero, so $\cap\,I_Z\subset \pi A$.
As the quotients $A/I_Z$ have no $\pi$-torsion, it follows that
$\cap\, I_Z\subset \pi^nA$ for all $n$. As $A$ is separated for the $\pi$-adic
topology, we get $\cap\, I_Z=0$. So we have a diagram with all morphisms injective
$$
\xymatrix@R=16pt{\hat{A}\; \ar[r] & \hat{A}_K \\
A\; \ar[r] \ar[u] & \ar[u] A_K
}
$$
Moreover we have $\hat{A}\cap A_K=A$. Indeed an element in $\hat{A}\cap A_K$ is
written $a/\pi^d$ with $a\in A$ and $d$ minimal, such that there exists a
system $(a_Z)_{Z\in \clF}$ with $a\equiv \pi^d a_Z$ modulo $I_Z$ for all
$Z$. If $d\geq 1$, reducing modulo $\pi$ we get $a\equiv 0$ modulo
$I_Z\otimes k$. From $\cap\, I_Z\otimes k=0$ it
follows that $a\in \pi A$, which is a contradiction. Hence $a/\pi^d\in A$.

We write $\mu_Z^\sharp\colon A/I_Z\to R\clG\otimes A/I_Z$ for the
coactions (step (2)). As $R\clG$ is finite free over $R$,
tensor product and inverse limit commute so we can pass to the limit and
obtain
$$
\hat{\mu}^\sharp:=\underset{\longleftarrow}{\lim}\,\mu_Z^\sharp
\colon \hat{A}\to R\clG\otimes\hat{A}
$$
On the generic fibre, all actions are induced from $\mu$ so that
$\hat{\mu}^\sharp\otimes K$ maps $A_K$ into $R\clG\otimes A_K$. From
$\hat{A}\cap A_K=A$ it follows that $\hat{\mu}^\sharp$ maps $A$ into
$R\clG\otimes A$. Thus $\hat{\mu}^\sharp$ extends to an action
of $\clG$ on $X$, universally faithful by construction.
\end{proo}

For algebraic schemes, we will use theorem \ref{theoZ} in the following
more convenient version.

\begin{coro} \label{coroMT}
Let $G$ be a finite flat group scheme over $R$. Let $X$ be a flat scheme of
finite type over $R$ and let $\mu\colon G\times X\to X$ be an action. We
assume that $X$ is covered by $G$-stable open affines $U_i$ with function
ring separated for the $\pi$-adic topology, such that $G$ acts faithfully
on the generic fibre $U_{i,K}$.
Then, if $X$ has reduced special fibre, there exists an effective model
for the action of $G$.
\end{coro}

\begin{proo}
If we have an effective model $\clG^h$ after extension to the
henselization $R^h$, then by unicity it descends to $R$ as well as the
action. So we may assume that $R$ is henselian. We will prove that
the family $\clF_k$ (notation of \ref{theoZ}) contains the set of
Cohen-Macaulay points, which is dense open. Since $X_k$ is reduced, it
follows that $\clF_k$ is schematically dense and we can apply theorem
\ref{theoZ}.
So let $x\in X_k$ be Cohen-Macaulay and let $(r_i)_{1\leq i\leq m}$ be
a system of parameters for the ring $\clO_{X_k,x}$. On an affine
neighbourhood $U=\Spec(A)$ of $x$ in $X$,
pick sections $s_1,\dots,s_m\in A$ whose germ at $x$ map to $r_i$ in
$\clO_{X_k,x}$. Let $Y$ be the closed subscheme of $U$ defined by the
vanishing of the $s_i$. As $(r_i)$ is a regular sequence, it follows
that $\clO_{Y,x}$ is flat over $R$.
Furthermore, $\clO_{Y_k,x}$ is artinian (by the Cohen-Macaulay
assumption) hence $\clO_{Y,x}$ is also quasi-finite over $R$.
Since $R$ is henselian $\clO_{Y,x}$ is even finite over $R$. Hence
the schematic image of $\Spec(\clO_{Y,x})$ in $X$ is a subscheme
$Z\in \clF$ which contains $x$.
\end{proo}

\subsection{Applications to Raynaud's group scheme}

Here we explain the link between our effective model and the
constructions in \cite{Ab}, in particular we derive some
new features of the so-called Raynaud's group scheme presented
therein.

\begin{cons}
Let $X$ and $G$ be as in corollary \ref{coroMT}. Let $f\colon X\to
Y=X/G$ be the quotient morphism. Let $Z_1,\dots,Z_r$ be the orbits of
irreducible components of $X_k$ under the action of $G$. Let $U_i$
denote the open subscheme of $X$ obtained by removing all the
components of $X_k$ except those in the orbit $Z_i$. Let
$$
X^*:=\cup \, U_i
$$
In $Y$ we consider $V_i:=f(U_i)$ and $Y^*:=f(X^*)$.
By theorem \ref{coroMT} there is a finite flat group scheme $\clG_i\to
\Spec(R)$ which is an effective model for the action of $G$ on $U_i$.
Let $\clG$ be the flat group scheme over $Y^*$ obtained by glueing the
schemes $\clG_i\times V_i$ along $X_K$. We call it {\em Raynaud's group
scheme}. \qedend
\end{cons}

Note that $X^*$ and $Y^*$ strictly contain the smooth loci of
$X$ and $Y$~: actually the only singularities that are omitted are
the intersections of components of the special fibre. Moreover, under
the assumptions of theorem 3.1.1 of \cite{Ab}, by unicity it is
clear that our group scheme $\clG\to Y^*$ coincides with Raynaud's
group scheme as defined in \cite{Ab} over the smooth
locus. Thus we obtain~:

\begin{coro} \label{precisions_on_RGS}
Let $X$ be a stable curve over $R$ with smooth generic fiber, $G$ a
finite group acting on $X$, and $Y=X/G$. Assume that the closure of
fixed points of $G$ in $X_K$ are disjoint sections lying in the smooth
locus $X_\sm$. Assume $p^2{\not |}\: |G|$ and the $p$-Sylow subgroup
of $G$ is normal. Let $\clG\to Y_\sm$ be Raynaud's group scheme as in
theorem 3.1.1 of {\rm \cite{Ab}} (notations $X$ and $Y$ are inverted
there). Then $\clG$ is constant on (orbits of) irreducible components
of the special fibre, and it extends to nodes $p\in Y_k$ lying on one
single component. \qedend
\end{coro}

We stress that the fact that $\clG$ extends to some nodes does not
mean that the extension of $\clG$ defined in \cite{Ab}, \S~3.2 is
representable by a scheme at these points, because the latter is
actually endowed with a supplementary stack structure so as to make
the covering into a torsor (always assuming that $p^2$ does not divide
$|G|$).

When $|G|$ is arbitrary, $X^*\to Y^*$ will not be a torsor under
$\clG$, because general points in the special fibre of $X^*$ may have
nontrivial stabilizers. Note that this prevents any hope of
commutation between quotient and base change. We will see an example
of this for $G=\zmod{p^2}$ in the next section.

\section{Examples}

Computations of effective models provide a whole zoo of
finite flat group schemes. Effective models of $(\zmod{p})_R$ appear
in \cite{OSS}, I. \S~2 or \cite{He}, \S~1 for degenerations of
$\mu_p$-torsors in unequal characteristics and in \cite{Ma}, \S~3.2 for
degenerations of $\zmod{p}$-torsors in equal characteristic $p>0$.

Now we give examples in degree $p^2$. Recently Mohamed Sa{\"\i}di
studied degenerations of torsors under $\zmod{p^2}$ in equal
characteristics \cite{Sa2}. He computed equations for such degenerations~; here
we will provide the group scheme which extends the $\zmod{p^2}$ action
on this set of equations. We will study one case where one gets a torsor
structure, and one where this fails to happen (of course all cases in
theorem 2.4.3. of \cite{Sa2} could be treated similarly).

Note that it is very likely that one could give similar examples in the case
of mixed characteristics, using the Kummer-to-Artin-Schreier isogeny of
Sekiguchi and Suwa in degree $p^2$. The computations would just be (possibly
substantially) more complicated.

\subsection{Witt vectors of length 2}

We assume that $R$ is complete and has equal characteristics $p>0$,
so $R\simeq k[[\pi]]$. Under this assumption, torsors under
$\zmod{p^2}$ are described by Witt theory.

\begin{noth}
{\bf Classical Witt theory.}
First we briefly recall the notations of Witt theory in degree $p^2$
(see \cite{DG}, chap. V). The group scheme of Witt vectors of length 2
over $R$ has underlying scheme
$W_{2,R}=\Spec(R[u_1,u_2])\simeq \bbA^2_R$ with multiplication law
$$
(u_1,u_2)+(v_1,v_2)=\big(u_1+v_1,u_2+v_2+\sum_{k=1}^{p-1} \bin{p}{k} \, u_1^kv_1^{p-k}\big)
$$
Here we put once for all $\bin{p}{k}:=\frac{1}{p}{p\choose k}$ where
${p\choose k}$ is the binomial coefficient.
The Frobenius morphism of $W_2$ is denoted by $F(u_1,u_2)=(u_1^p,u_2^p)$.
Put $\phi:=F-\Id$. From the exact sequence
$$
0\to (\zmod{p^2})_R\to W_{2,R}\stackrel{\phi}{\longrightarrow} W_{2,R}\to 0
$$
it follows that any {\'e}tale torsor $f\colon \Spec(B)\to \Spec(A)$
under $(\zmod{p^2})_R$ is given by an equation
$$
F(X_1,X_2)-(X_1,X_2)=(a_1,a_2)
$$
where $(a_1,a_2)\in W_2(A)$ is a Witt vector and the substraction
is that of Witt vectors. Furthermore, $(a_1,a_2)$ is well-defined up
to addition of elements of the form $F(c_1,c_2)-(c_1,c_2)$. Note that
$$
F(X_1,X_2)-(X_1,X_2)=\big(X_1^p-X_1,X_2^p-X_2+
\sum_{k=1}^{p-1} \bin{p}{k} \, (X_1)^{pk}(-X_1)^{p-k}\big)
$$
We emphasize that the Hopf algebra of $(\zmod{p^2})_R$ is
$$
R[\zmod{p^2}]=\frac{R[u_1,u_2]}{(u_1^p-u_1,u_2^p-u_2)}
$$
with comultiplication that of $W_2$.
\end{noth}

\begin{noth} \label{twistedforms}
{\bf Twisted forms of $W_2$.}
Let $\lambda,\mu,\nu$ be elements of $R$. We define a "twisted" group
$W_2^\lambda$ as the group with underlying scheme $\Spec(R[u_1,u_2])$ and
multiplication law given by
$$
(u_1,u_2)+(v_1,v_2)=\bigg(u_1+v_1\, ,\, u_2+v_2+\lambda\,
\sum_{k=1}^{p-1} \bin{p}{k} \, u_1^kv_1^{p-k}\bigg)
$$
We have the following analogues of the scalar multiplication and the
Frobenius of $W_2$ :
$$
\begin{array}{rl}
I_{\lambda,\mu}^{\nu} \, \colon \; W_2^\lambda & \longrightarrow W_2^{\lambda\mu}
\medskip \\
(u_1,u_2) & \longmapsto (\nu u_1,\mu\nu^p u_2) \\
\end{array}
$$
and
$$
\begin{array}{rl}
F_{\lambda} \, \colon \; W_2^\lambda & \longrightarrow W_2^{\lambda^p}
\medskip \\
(u_1,u_2) & \longmapsto (u_1^p,u_2^p) \\
\end{array}
$$
In case $\mu=\lambda^{p-1}$ we define an isogeny
$$
\phi_{\lambda,\nu}:=F_{\lambda}-I_{\lambda,\lambda^{p-1}}^{\nu}\colon
W_2^\lambda \to W_2^{\lambda^p}$$
We have
$$
\phi_{\lambda,\nu}(u_1,u_2) =
\bigg(u_1^p-\nu u_1\, ,\, u_2^p-\nu^p\lambda^{p-1}u_2
+\lambda^p \sum_{k=1}^{p-1} \bin{p}{k} \, u_1^{pk}(-\nu u_1)^{p-k} \bigg)
$$
The kernel $\clK_{\lambda,\nu}:=\ker(\phi_{\lambda,\nu})$
is a finite flat group of rank $p^2$. If $p>2$ its Hopf algebra is
$$
R[\clK_{\lambda,\nu}]=\frac{R[u_1,u_2]}{(u_1^p-\nu u_1,u_2^p-\nu^p\lambda^{p-1}u_2)}
$$
\end{noth}

\subsection{Two examples}

The examples come from the following situation studied by Sa{\"\i}di.
Denote by $G=\zmod{p^2}$ the constant group, and by $Y=\bbA^1_R=\Spec(R[w])$
the affine line over $R$. Let $m_1,m_2\in \bbZ$ be integers. Let
$f_K\colon X_K\to Y_K$ be the $(\zmod{p^2})_K$-torsor
over $Y_K=\bbA^1_K$ given by the equations :
$$
\left\{
\begin{array}{lcl}
T_1^p-T_1 & = & \pi^{m_1}w \\
T_2^p-T_2 & = & \pi^{m_2}w - {\displaystyle \sum_{k=1}^{p-1}} \bin{p}{k} \, (T_1)^{pk}(-T_1)^{p-k} \\
\end{array}
\right.
$$
Depending on the values of the conductors
$m_1$, $m_2$ this gives rise to different group degenerations.

\begin{exam}
Assume $m_1=0$ and $m_2=-p$. Then after the change of variables $Z_1=T_1$,
$Z_2=\pi T_2$ the map $f_K$ extends to a cover $X\to Y$ with equations
$$
\left\{
\begin{array}{lcl}
Z_1^p-Z_1 & = & w \\
Z_2^p-\pi^{(p-1)}Z_2 & = & w - \pi^p{\displaystyle \sum_{k=1}^{p-1}} \bin{p}{k} \, (Z_1)^{pk}(-Z_1)^{p-k} \\
\end{array}
\right.
$$
It is quickly seen that the action of $\zmod{p^2}$ extends to $X$.
As is obvious from the expression of the isogeny $\phi_{\lambda,\nu}$
(see \ref{twistedforms}), $X\to Y$ is a torsor under $\clK_{\lambda,\nu}$
for $\lambda=\pi$ and $\nu=1$. Thus, the effective model is $G_R\to
\clG$ with $\clG=\clK_{\pi,1}$.
\end{exam}

\begin{exam}
Assume $m_1=-p^2n_1<0$ and $m_2=0$. Put $\tilde{m}_1=n_1(p(p-1)+1)$.
Then after the change of variables $Z_1=\pi^{pn_1} T_1$ and
$Z_2=\pi^{\tilde{m}_1}T_2$ the map $f_K$ extends to a cover $X\to Y$
with equations
$$
\left\{
\begin{array}{lcl}
Z_1^p-\pi^{(p-1)pn_1}Z_1 & = & w \\
Z_2^p-\pi^{(p-1)\tilde{m}_1}Z_2 & = & \pi^{p\tilde{m}_1}w
- {\displaystyle \sum_{k=1}^{p-1}} \bin{p}{k} \, \pi^{pn_1(p-1)(p-1-k)}(Z_1)^{pk}(-Z_1)^{p-k} \\
\end{array}
\right.
$$
The action of $\zmod{p^2}$ extends to this model as follows :
for $(u_1,u_2)$ a point of $G_R=(\zmod{p^2})_R$,
$$
(u_1,u_2).(Z_1,Z_2)=
\bigg(Z_1+\pi^{pn_1}u_1\; ,\; Z_2+\pi^{\tilde{m}_1}u_2+
\sum \bin{p}{k} \, \pi^{n_1(p(p-1)+1-pk)}(Z_1)^k(u_1)^{p-k}\bigg)
$$
In order to find out the model $G_R\to \clG$ we look at the subalgebra of $RG$
generated by $v_1=\pi^{n_1}u_1$ and $v_2=\pi^{\tilde{m}_1}u_2$ :
$$
R\clG:=R[v_1,v_2]\subset RG
$$
One computes that $R\clG$ inherits a comultiplication from $RG$ :
$$
(v_1,v_2)+(w_1,w_2)=
\bigg(v_1+w_1\; ,\; v_2+w_2+ \sum \bin{p}{k} \, \pi^{n_1(p-1)^2} v_1^kw_1^{p-k}\bigg)
$$
Thus if $p>2$ we recognize $\clG\simeq \clK_{\lambda,\nu}$ for
$\lambda=\pi^{n_1(p-1)^2}$ and $\nu=\pi^{n_1(p-1)}$. The action
of $G$ on $X$ extends to an action of $\clG$ as
$$
(v_1,v_2).(Z_1,Z_2)=
\bigg(Z_1+\pi^{(p-1)n_1}v_1\; ,\; Z_2+v_2+
\sum \bin{p}{k} \, \pi^{n_1(p-1)(p-1-k)}Z_1^kv_1^{p-k}\bigg)
$$
Here $X\to Y$ is not a torsor under $\clG$. Indeed, on the special fibre
we have $\clG_k=(\alpha_p)^2$ and the action on $X_k$ is
$$
(v_1,v_2).(Z_1,Z_2)=
\bigg(Z_1\; ,\; Z_2+v_2+v_1Z_1^{p-1}\bigg)
$$
This action is faithful as required, but any point $(z_1,z_2)\in X_k$ has a
stabilizer of order $p$ which is the subgroup of $\clG_k$ defined by the equation
$v_2+v_1 z_1^{p-1}=0$.
\end{exam}

\v{5}

\begin{flushleft}
{\em
Institut de Math{\'e}matiques

Analyse Alg{\'e}brique

Universit{\'e} Pierre et Marie Curie

Case 82

4, place Jussieu

F-75252 Paris Cedex 05

romagny@math.jussieu.fr}
\end{flushleft}

\end{document}